\documentclass[12pt]{article}

\usepackage{amsmath}
\usepackage{amssymb}
\usepackage{amsthm}
\usepackage[dvips]{graphicx}
\usepackage{color}
\usepackage{epsfig}
\usepackage[mathscr]{eucal}
\newtheorem{theorem}{Theorem}

\newtheorem{prop}[theorem]{Proposition}
\newtheorem{conj}{Conjecture}
\newtheorem{lemma}[theorem]{Lemma}
\newtheorem{cor}[theorem]{Corollary}

\newtheorem{definition}[theorem]{Definition}

\allowdisplaybreaks

\def\ispreprint{}

\ifx\ispreprint\undefined

\else
\fi
\title{Quasirandom Arithmetic Permutations}
\author{Joshua N. Cooper\\ \small Department of Mathematics \\ \small Courant Institute of Mathematical Sciences \\ \small New York University, New York, NY \\ \small \texttt{cooper@cims.nyu.edu}}
\date{\today}

\begin{document}

\maketitle

\ifx\ispreprint\undefined
Correspondence Information: \\
Dr. Joshua N. Cooper, 312 2nd Ave West, Apt 414, Seattle, WA 98119 \\
Email: cooper@cims.nyu.edu \\
Phone: (858) 699-6049
\pagebreak
\else
\fi

\begin{abstract}
In \cite{C1}, the author introduced {\it quasirandom permutations}, permutations of $\mathbb{Z}_n$ which map intervals to sets with low discrepancy.  Here we show that several natural number-theoretic permutations are quasirandom, some very strongly so.  Quasirandomness is established via discrete Fourier analysis and the Erd\H{o}s-Tur\'{a}n inequality, as well as by other means.  We apply our results on S\'{o}s permutations to make progress on a number of questions relating to the sequence of fractional parts of multiples of an irrational.  Several intriguing open problems are presented throughout the discussion.
\end{abstract}

\ifx\ispreprint\undefined
Keywords: \\
Quasirandom, permutation, discrepancy, exponential sums, S\'{o}s permutation.
\pagebreak
\else
\fi
\section{Introduction}

Random objects play a crucial role in modern combinatorial theory.  Despite their success, it is often desirable to replace random constructions with explicit ones.  One approach to doing so involves quantifying the randomness of a given object in a manner which is relevant to the structure under consideration, and then constructing objects which resemble random ones in this metric.  One particularly useful measure of randomness is discrepancy, which measures how uniformly an object's substructures are organized.  If, in the limit, an infinite sequence of objects has high uniformity, i.e., low discrepancy, we call it ``quasirandom''.

Presently, we apply a discrepancy-theoretic definition of quasirandomness to several types of very natural arithmetic permutations.  Riding on the tails of extensive work concerning equidistribution properties of exponential sums (e.g., \cite{K1,KS1,LN1,N1}) and the uniform distribution of the sequence $\{n \alpha\}$, for $\alpha$ irrational (e.g, \cite{DT1,O1,Sg1}), this analysis provides additional justification for the intuition that these functions are all random-like.

The maps under consideration send $\mathbb{Z}_n$ (or $\mathbb{Z}_n^\times$) to itself as follows.  We assume that $p$ is prime.

\begin{enumerate}
\item $\psi_k$ : For $k \in \mathbb{Z}_n^\times$ and $s \in \mathbb{Z}_n$, define $\psi_k(s) = ks$.
\item $\lambda_a$ : For $s \in \mathbb{Z}_p^\times$ and $a \in \mathbb{Z}_p^\times$, $\lambda_a(s) = as^{-1}$.
\item $\eta_{a,k}$ : For $s \in \mathbb{Z}_p$, $a \in \mathbb{Z}_p^\times$, and $k \in \mathbb{Z}_{p-1}$ with $(k,p-1)=1$, $\eta_k(s) = as^k$.
\item $\rho_{a,\tau}$ : For $\tau$ a primitive root of $\mathbb{Z}_p$ and $a \in \mathbb{Z}_p^\times$, $\rho_{a,\tau}(s) = a \, \tau^s$.
\item $\beta_\alpha$ : For $s,t \in \{1,\ldots,n\}$ and $\alpha \in \mathbb{R}$ irrational, $\beta_\alpha(s) < \beta_\alpha(t)$ iff $\{\alpha s\} < \{\alpha t\}$, where $\{x\}$ is the fractional part of $x$.
\end{enumerate}

The $\beta_\alpha$ we will refer to as S\'{o}s permutations, in honor of V. S\'{o}s' early and definitive work on them, as typified by \cite{So1} and \cite{So2}.  Throughout this chapter, we will notationally suppress the dependence of the functions defined above on $n$ and $p$.

The author defined quasirandom permutations in \cite{C1}, following Chung, Graham, and Wilson's introduction of quasirandom graphs, hypergraphs, tournaments, subsets of $\mathbb{Z}_n$, and others (\cite{CG1,CG91ss,CG3,CG92sz,CGW1}).  The central thrust of that paper, as with other studies of quasirandomness, was the demonstration that several different measures of randomness are really the same, in that being random in any one of these ways is equivalent to all the others.  Thus, establishing quasirandomness requires only that one of these properties be shown to hold.

In the next section, we give formal definitions and discuss the results of \cite{C1} briefly.  We also provide a tool which will be needed later, the Erd\H{o}s-Tur\'{a}n inequality.  In Section \ref{proofs}, we show that the above permutations are quasirandom for broad ranges of parameters.  Each case is accompanied by open questions concerning the true order of magnitude of their discrepancies.  We conclude in Section \ref{appsfrac} with an application to a question of K. O'Bryant concerning the sequence of fractional parts of integer multiples of an irrational.

\section{Preliminaries} \label{prelims}

We consider permutations, elements of $S_n$, to be actions on $\mathbb{Z}_n$ or on the integers $[n]=\{1,\dots,n\}$, and we will write them in ``one-line'' notation when needed, i.e., $\sigma \in S_n$ will be written $(\sigma(0) \, \ldots \, \sigma(n-1))$.  An \textit{interval} of $\mathbb{Z}_n$ is any set $I$ which is the projection of an interval of $\mathbb{Z}$.  (Note that this definition permits ``wrap-around''.)  For any $S, T \subset \mathbb{Z}_n$ we define the {\it discrepancy} of $S$ in $T$ as
$$
D_T(S) = \left | | S \cap T | - \frac{|S||T|}{n} \right |,
$$
and we define the discrepancy of a permutation $\sigma$ by
$$
D(\sigma) = \max_{I,J} D_J(\sigma(I)),
$$
where $I$ and $J$ vary over all intervals of $\mathbb{Z}_n$.  Also, define
$$
D^*(\sigma) = \max_{I,J} D_J(\sigma(I)),
$$
where $I$ and $J$ vary only over ``initial'' intervals, i.e., projections of intervals of the form $[0,M]$ for $M \geq 0$.  Then, since discrepancy is ``subadditive'', it is clear that $D^*(\sigma) \leq D(\sigma) \leq 4 D^*(\sigma)$ for any $\sigma$.

We say that a sequence $\{\sigma_i\}_{i=1}^\infty$ of permutations of $\mathbb{Z}_{n_1}$, $\mathbb{Z}_{n_2}$, $\ldots$ is {\it quasirandom} if $D(\sigma_i) = o(n_i)$.  Furthermore, we will often suppress the indices and simply say that $D(\sigma) = o(n)$.  Thus, it is easy to see that $\sigma$ is quasirandom iff $\sigma^{-1}$ is.

Define ${\bf X}^\tau(\sigma)$ for $\tau \in S_m$ and $\sigma \in S_n$ to be the number of ``occurrences'' of $\tau$ in $\sigma$, i.e., the number of subsets $\{x_1 < \ldots < x_m\} \subset \mathbb{Z}_n$ such that $\sigma(x_i) < \sigma(x_j)$ iff $\tau(i) < \tau(j)$.  Throughout the rest of this chapter, by $e(x)$, we mean $e^{2 \pi i x}$, and by $f(n) \ll g(n)$, we mean that there exists a $C$ so that, for sufficiently large $n$, $f(n) \leq C g(n)$.  The following theorem appears in \cite{C1}.  We use the convention that the name of a set and its characteristic function are the same.

\begin{theorem} For any sequence of permutations $\sigma \in S_n$, integer $m \geq 2$, and fixed real $\alpha > 0$, the following are equivalent:
\flushleft{\begin{tabular}{ll}
{\bf [UB]} & \parbox[t]{5.0in}{(Uniform Balance) $D(\sigma) = o(n)$.} \\[.13in]
{\bf [SP]} & \parbox[t]{5.0in}{(Separability) For any intervals $I,J,K,K^\prime \subset
\mathbb{Z}_n$,
$$
\left | \sum_{x \in K \cap \sigma^{-1}(K^\prime)} I(x)
J(\sigma(x)) - \frac{1}{n} \sum_{x \in K,y \in K^\prime} I(x)J(y)
\right | = o(n)
$$} \\[.13in]
{\bf [mS]} & \parbox[t]{5.0in}{(m-Subsequences) For any permutation $\tau \in S_m$ and intervals
$I,J \subset \mathbb{Z}_n$ with $|I| \geq n/2$ and $|J| \geq n/2$, we have $|I \cap \sigma^{-1}(J)| \geq n/4 + o(n)$ and
$$
{\bf X}^\tau(\sigma |_{I \cap \sigma^{-1}(J)}) = \frac{1}{m!}
\binom{|\sigma(I) \cap J|}{m} + o(n^m).
$$} \\[.13in]
{\bf [2S]} & \parbox[t]{5.0in}{(2-Subsequences) For any intervals $I,J \subset \mathbb{Z}_n$ with $|I| \geq n/2$ and $|J| \geq n/2$, we have $|I \cap \sigma^{-1}(J)| \geq n/4 + o(n)$ and
$$
{\bf X}^{(01)}(\sigma |_{I \cap \sigma^{-1}(J)}) - {\bf
X}^{(10)}(\sigma |_{I \cap \sigma^{-1}(J)}) = o(n^2).
$$} \\[.13in]
{\bf [E($\alpha$)]} & \parbox[t]{5.0in}{(Eigenvalue Bound $\alpha$) For all nonzero $k \in \mathbb{Z}_n$ and any interval $I$,
$$
\sum_{s \in \sigma(I)} e(-ks/n) = o(n|k|^{\alpha}).
$$} \\[.13in]
{\bf [T]} & \parbox[t]{5.0in}{(Translation) For any intervals $I,J$,
$$
\sum_{k \in \mathbb{Z}_n} \left ( |\sigma(I) \cap (J+k)| - \frac{|I||J|}{n} \right )^2 = o(n^3).
$$}
\end{tabular}}

Furthermore, for any implication between a pair of properties above, there exists a constant $K$ so that the error term $\epsilon_2 n^k$ of the consequent is bounded by the error term $\epsilon_1 n^l$ of the antecedent in the sense that $\epsilon_2 \ll \epsilon_1^K$.
\end{theorem}

This result is interesting particularly because it says that once we show one of these properties for a sequence of permutations, we get the rest for free.  For example, we will show that the permutation $\lambda_1$ has discrepancy at most $p^{1/2+\epsilon}$ (i.e, property ${\bf [UB]}$), so it also has approximately the ``right number'' of inversions, i.e., ${\bf [2S]}$.  In other words, $x^{-1} < x$ about as often as $x^{-1} > x$.  Furthermore, using the quantitative statement of the theorem, we can show that the difference between the numbers of $x$'s satisfying these two conditions is $\ll p^{1/2+\epsilon}$.  It is also a simple matter to show the following:

\begin{prop} Define $A = \sqrt{n D(\sigma)}$.  If $I$ and $J$ are two intervals of $\mathbb{Z}_n$ with $|I|>A$ and $|J|>A$, then $\sigma(I) \cap J \neq \emptyset$.
\end{prop}

It follows immediately that, for any $\epsilon > 0$ and sufficiently large $p$, if we have two intervals of length at least $p^{3/4+\epsilon}$, there is a point in one whose inverse mod $p$ lies in the other.

The following result giving a universal lower bound on the discrepancy of a permutation appears in \cite{C1} and follows immediately from a result of W. Schmidt \cite{S1}.

\begin{prop} \label{disclowerbound2} For $\sigma \in S_n$, $D(\sigma) = \Omega(\log n)$.
\end{prop}

This bound is actually achievable, since if $\sigma$ is taken to be the permutation which reverses the binary expansion of integers between $0$ and $2^n-1$, then $D(\sigma) \ll n$ (q.v. \cite{C1}).  Therefore there are ``maximally'' quasirandom sequences of permutations, and then other ones whose discrepancies grow faster than $\log n$.  Interestingly, random permutations can be shown to have discrepancy $\ll \sqrt{n \log n}$, and it is straightforward to show that for almost all permutations $\sigma$, $D(\sigma) \gg \sqrt{n}$.  The phenomenon of random objects being less uniform than specially constructed ones is a common phenomenon, and appears throughout combinatorics, discrepancy theory, the theory of quasi-Monte Carlo integration, and elsewhere.

We have the following standard lemma, which will be needed later.

\begin{lemma} \label{intb2} If $J$ is an interval of $\mathbb{Z}_n$, then $\tilde{J}(k) \leq \frac{n}{2|k|}$.
\end{lemma}
\begin{proof}
We may write the magnitude of the $k^{\mbox{\scriptsize th}}$
Fourier coefficient of $J = [a+1,a+M]$ as
\begin{align*}
|\tilde{J}(k)| &= | \sum_x J(x) e(-kx/n) | = | \sum_{x=a}^b e(-kx/n) | = | \sum_{x=1}^M e(-kx/n) | \\
& = \frac{|e(-kM/n) - 1|}{|e(-k/n) - 1|} \leq \frac{2}{4 |k| / n} = \frac{n}{2|k|}
\end{align*}
since $|e^{i \theta} - 1| \geq \frac{2 |\theta|}{\pi}$ for all
$\theta$.
\end{proof}

Finally, we present the Erd\H{o}s-Tur\'{a}n inequality (\cite{ET1}), which gives a bound on the discrepancy of a sequence in terms of its Fourier transform.  The discrepancy of a sequence $S=\{x_i\}_{i=0}^{m-1}$ of reals in $[0,1)$, i.e., elements of $\mathbb{R}/\mathbb{Z}$, is defined to be
$$
D(S) = \sup_{0 \leq \alpha \leq 1} \big | |\{i : 0 \leq i < m, x_i \in [0,\alpha) \}| - \alpha m \big |.
$$

\begin{theorem}[Erd\H{o}s-Tur\'{a}n, 1948] \label{ETI} For a sequence $\{x_i\}_{i=0}^{m-1} \subset \mathbb{R}/\mathbb{Z}$, define
$$
A(k)=\sum_{i=0}^{m-1} e(k x_i).
$$
Then there is an absolute constant $C$ so that, for any positive integer $K$,
$$
D(U) \leq C \left ( \frac{m}{K} + \sum_{k=1}^K \frac{|A(k)|}{k} \right ).
$$
\end{theorem}

If we take $x_i = \sigma(i)/n$ for some $\sigma \in S_n$, then we have the following version of Theorem \ref{ETI}.

\begin{cor} \label{CETI1} Let $\sigma \in S_n$, $n > 1$, and suppose that for all $m$ and $k$,
$$
\left | \sum_{s=0}^{m-1} e(k \sigma(s)/n) \right | \leq \alpha(n).
$$
Then there is an absolute constant $C_0$ so that
$$
D(\sigma) \leq C_0\,\alpha(n) \log n.
$$
\end{cor}
\begin{proof} Take $K = \left \lceil m / \alpha(n) \right \rceil$, and simplify.
\end{proof}

It is well known (see \cite{K1}, Theorem 2) that
$$
\left | \sum_{s=j+1}^{j+M} e(f(s)) \right | \leq \max_{0 \leq a \leq n-1} \left | \sum_{s=0}^{n-1} e(f(s) + as/n) \right | (1 + \log n)
$$
so we have another useful corollary of Theorem \ref{ETI}:

\begin{cor} \label{CETI2} Let $\sigma \in S_n$, $n>1$, and suppose that for all $a$ and $k\neq 0$,
$$
\left | \sum_{s=0}^{n-1} e\left(\frac{k \sigma(s)+as}{n}\right) \right | \leq \beta(n).
$$
Then there is an absolute constant $C_1$ so that
$$
D(\sigma) \leq C_1\,\beta(n) \log^2 n.
$$
\end{cor}

\section{Proofs} \label{proofs}

In this section, we show that the five permutations listed in Section \ref{prelims} are (usually) quasirandom.

\subsection{Multiplication}

We begin with $\psi_k$, used by Alon \cite{A1} to derandomize a maximum-flow algorithm of Cheriyan and Hagerup.  Recall the definition of $\psi_k$:

\begin{definition} For $k \in \mathbb{Z}^\times_n$, write $\psi_k$ for the permutation which sends $s \in \mathbb{Z}_n$ to $s\!\cdot\!k$.
\end{definition}

The following theorem says that multiplication by \textit{some} units of $\mathbb{Z}_n$ comes fairly close to meeting the lower bound of Schmidt.

\begin{theorem} \label{phik} For each $n$,
$$
\mathbf{E} [D(\psi_k)] = O(\log^2 n \log\log n),
$$
where the expected value is taken over all $k \in \mathbb{Z}_{n^\times}$. 
\end{theorem}
\begin{proof} Clearly, it suffices to show that the expected value of $D(\psi_k)$ is $O(\log^2{n} \cdot \log \log n )$ when we choose a random $k$ uniformly from $\mathbb{Z}^\times_n$.  Note that $\sum_s e(sx/n) = n \cdot \chi(x=0)$, where the sum is over all elements of $\mathbb{Z}_n$.  Thus, for any intervals $I=[a,b],J=[c,d] \subset \mathbb{Z}_n$,
$$
\left | \psi_k(I) \cap J \right | = n^{-1} \sum_s \sum_{y \in I} \sum_{z \in J} e(s(yk-z)/n).
$$
Since the term with $s=0$ is just equal to $|I||J|$, it is easy to see that
$$
D(\psi_k) = \sup_{I,J} \left | n^{-1} \sum_{s \neq 0} \sum_{y \in I} \sum_{z \in J} e(s(yk-z)/n) \right |.
$$
and, summing the resulting geometric series,
\begin{align*}
D(\psi_k) &= \sup_{I,J} \left | n^{-1} \sum_{s \neq 0} \left ( \sum_{y \in I} e(syk/n) \right ) \left ( \sum_{z \in J} e(-sz/n) \right ) \right | \\
&= \sup_{I,J} \left | n^{-1} \sum_{s \neq 0} \frac{e(sk(b+1)/n)-e(ska/n)}{e(sk/n)-1} \cdot \frac{e(-s(d+1)/n)-e(-sc/n)}{e(-s/n)-1} \right | \\
&\leq n^{-1} \sup_{I,J} \sum_{s \neq 0} \left | \frac{e(sk|I|/n)-1}{e(sk/n)-1} \right | \cdot \left | \frac{e(-s|J|/n)-1}{e(-s/n)-1} \right | \\
& \leq n^{-1} \sum_{s \neq 0} \frac{n^2}{4|s||sk|} = \sum_{s \neq 0} \frac{n}{4|s||sk|} \, \, ,
\end{align*}
by the proof of Lemma \ref{intb2}.  Taking the expected value of this sum over all $k \in \mathbb{Z}^\times_n$, we have
$$
\mathbf{E}[D(\psi_k)] \leq \frac{1}{\phi(n)} \sum_{k \in \mathbb{Z}^\times_n} \sum_{s \neq 0} \frac{n}{4|s||sk|}.
$$
where we have divided by the Euler $\phi$-function.  It is well known (\cite{HW1}) that $\phi(n) = \Omega(n/\log \log n)$, so we may conclude that
\begin{align*}
\mathbf{E}[D(\psi_k)] & \leq O(\log \log n) \sum_{s \neq 0} \sum_{k \neq 0} \frac{1}{4|s||sk|} \\
& = O(\log \log n) \left ( \sum_{s \neq 0} \frac{1}{|s|} \right ) \left ( \sum_{k \neq 0} \frac{1}{|sk|} \right ) \\
& = O(\log^2 n \log \log n).
\end{align*}
Therefore, there is some $k$ so that $D_J(\psi_k(I)) = O(\log^2 n \log \log n)$.
\end{proof}

Note that the above argument allows us to drop the ``$\log \log n$'' if
$n=p$ is prime.  Therefore $\psi_k$ is {\it highly} quasirandom, for almost
all $k$.  Based on extensive computational evidence, we believe that the
true order of magnitude of $\mathbf{E}[D(\psi_k)]$ is, in fact, $\log^2 n$
for almost all $k$, but we are unable to prove this.  Furthermore, computer
evidence points even more strongly to the existence of a $k$ for each $n$
with $D(\psi_k) = O(\log n)$.  In connection with ``good'' lattice points
for generating well-distributed points in the unit square, Neiderreiter
(\cite{N2}) has previously made this conjecture.  The best known bounds are
given by Larcher (\cite{La1}), who has shown that a $k$ always exists so
that $D(\psi_k) = O(\log n (\log \log n)^2)$.  We comment further on this
conjecture in Section \ref{Irrational}.

\begin{conj} \label{weaklogn} For some $k \in \mathbb{Z}_p^\times$, with $p$ prime, $D(\psi_k) = O(\log p)$.
\end{conj}

We would even venture the following stronger statement:

\begin{conj} \label{stronglogn}$\lim_{p \rightarrow \infty} \min_k D(\psi_k)/\log p = 1/2$.
\end{conj}

By virtue of Proposition \ref{disclowerbound2}, this would mean that there always exists a $k$ so that $\psi_k$ is maximally quasirandom.  The constant $1/2$ is close to the best known for maximally quasirandom permutations -- a result of H. Faure (\cite{F2}) on generalized van der Corput sequences implies the existence of a sequences of permutations $\sigma$ with $D(\sigma)/\log p \rightarrow 23/(35\log 6) \approx .367$.

\subsection{Exponentiation and Inversion}

The cases of exponentiation and inversion are particularly easy to deal with.  Let $p$ be a prime.  Recall the definitions of $\rho_{a,\tau}$ and $\lambda_a$:

\begin{definition} For $\tau$ a primitive root of $\mathbb{Z}_p$ and $a \in \mathbb{Z}_p^\times$, $\rho_{a,\tau}(s) = a \tau^s$.  (We define $\rho_{a,\tau}(0) = 0$ for convenience.)
\end{definition}

\begin{definition} For $s \in \mathbb{Z}_p^\times$ and $a \in \mathbb{Z}_p^\times$, $\lambda_a(s) = as^{-1}$.  (Again, let $\lambda_k(0)=0$.)
\end{definition}

The following theorem is usually known as the P\'olya-Vinogradov inequality:

\begin{theorem} For $\tau$ a primitive root of $\mathbb{Z}_p$,
$$
\left | \sum_{s=1}^{m} e(ka\tau^s/p) \right | \ll p^{1/2} \log p.
$$
uniformly in $m$ and $k$.
\end{theorem}

We may therefore conclude immediately, based on Corollary \ref{CETI1}, that

\begin{theorem} $D(\rho_{a,\tau}) \ll p^{1/2} \log^2 n$.
\end{theorem}

In this case, we believe the bound to be best possible (except possibly for the $\log$ terms).  Having taken care of ``exponentiation'' permutations, we can address ``inversion'' similarly.  The following classical result on Kloosterman sums, known as the Weil bound, appears in \cite{LN1}.  Define
$$
K(a,b)=\sum_{s \in \mathbb{Z}_p^\times} e((as + bs^{-1})/p).
$$

\begin{theorem} $|K(a,b)| \leq 2p^{1/2}$ if $b$ is nonzero.
\end{theorem}

We therefore have

\begin{theorem} $D(\lambda_k) \ll p^{1/2} log^2 p$.
\end{theorem}
\begin{proof} Note that, if $s=0$, $e((as + \lambda_k(s))/p) = 1$, so
$$
\left | \sum_{s \in \mathbb{Z}_p} e((as + \lambda_k(s))/p) \right | \leq 2p^{1/2} + 1
$$
and the result follows from Corollary \ref{CETI2}.
\end{proof}

Again, we conjecture that this bound is best possible, up to a possible $\log$ power.

\subsection{Powers}

Recall the definition of $\eta_{a,k}$:

\begin{definition} For $s \in \mathbb{Z}_p$, $a \in \mathbb{Z}_p^\times$, and $k \in \mathbb{Z}_{p-1}^\times$, define $\eta_{a,k}(s) = as^k$.
\end{definition}

It is an old and well known result of A. Weil that
$$
\left | \sum_{s=1}^p e(f(x)/p) \right | \leq (deg(f) - 1) \sqrt{p}
$$
for any $f \in \mathbb{Z}_p[x]$ and prime $p$.  This bound can be strengthened, however, if we restrict our attention to certain types of polynomials.  For example, the following result appears in \cite{Ka1}:

\begin{theorem} \label{kara} Let $f(x) = ax^k + bx \in \mathbb{Z}_p[x]$ with $a,b$ nonzero and $2 \leq k \leq p-1$.  Then
$$
\left | \sum_{s=1}^p e(f(x)/p) \right | \leq (n-1)^{1/4} p^{3/4}.
$$
\end{theorem}

Therefore, we have

\begin{theorem} If $p$ is a prime, $a \in \mathbb{Z}_p^\times$, and $k \in \mathbb{Z}_{p-1}$ with $(k,p-1)=1$ and $k \geq 2$, then
$$
D(\eta_{a,k}) \ll k^{1/4} p^{3/4} \log^2 p.
$$
uniformly in $a$.
\end{theorem}
\begin{proof} If $(k,p-1)=1$, then clearly $| \sum_{s=1}^p e(a x^k/p) | = 0$.  The result then follows from Theorem \ref{kara} by applying Corollary \ref{CETI2}.
\end{proof}

We can immediately deduce the following.

\begin{cor} Suppose $p$ is a prime, $a \in \mathbb{Z}_p^\times$, and $k \in \mathbb{Z}_{p-1}$ with $(k,p-1)=1$ and $k \geq 2$.  If $k = o(p/\log^8 p)$, then $\eta_{a,k}$ is quasirandom.
\end{cor}

In another direction, we can show that almost all exponents $k$ yield quasirandom permutations.  Define
$$
W_{a,c}(t) = \sum_{k=1}^{t} \left | \sum_{x=1}^{t} e((a \vartheta^x + c \vartheta^{xk})/p) \right |
$$
where $\vartheta$ is an integer of multiplicative order $t$ in $\mathbb{Z}_p$, $p \geq 3$ prime.  Also, let $d(k)$ denote the number of divisors of $k$.  The following theorem is proved in \cite{CFKLLS}.

\begin{theorem} For any $a,c \in \mathbb{Z}_p$, $c \neq 0$,
$$
W_{a,c}(t) \ll \left \{ \begin{array}{ll} t p^{1/2} d(t)\text{,} & \textrm{if } a = 0\text{;} \\ t^{5/3}p^{1/4}\text{,} & \text{otherwise.} \end{array} \right .
$$
\end{theorem}

Choose $\vartheta$ to be a primitive root, so that $\vartheta^x$ varies over all nonzero elements of $\mathbb{Z}_p$.  Using the fact that $d(k) \ll k^\epsilon$ for any $\epsilon > 0$, we may conclude that
$$
\sum_{k=1}^{p-1} \left | \sum_{x=1}^{t} e((ax + cx^k)/p) \right | \ll p^{23/12}.
$$
Therefore, if we choose $k$ randomly and uniformly from the $\phi(p-1)$ elements of $\mathbb{Z}_p^\times$, the expected size of $|\sum_{x=0}^{p-1} e((ax + cx^k)/p) - 1|$ is $O(p^{11/12} \log \log p)$, uniformly in $c \neq 0$.  Formally,

\begin{theorem} \label{almostpower} For almost all $k \in \mathbb{Z}_p^\times$ with $(k,p-1)=1$ and $k \geq 2$,
$$
D(\eta_{a,k}) \ll p^{11/12} \log \log p.
$$
uniformly in $a$ and $k$.
\end{theorem}

One might ask whether $S(a,k,M) = \sum_{s=1}^M e(a s^k/p) \ll p^{1-\epsilon}$ for some $\epsilon > 0$ uniformly in $a$, $k$, and $M$, since this would imply by Corollary \ref{CETI1} that $D(\eta_{a,k})$ is {\it always} quasirandom for $(k,p-1)=1$ and $k \geq 2$.  The answer is, unfortunately, no.  A result of Karacuba \cite{Ka1} states that, for some $k$ in the vicinity of $p/\log p$ there exists an $a \neq 0$ so that the Gauss sum $S(a,k,p) = p(1-o(1))$.  It is intriguing, however, that when $a$ and $b$ are nonzero, it is known that $|\sum_{s=1}^p e((as^k + bs)/p)| \leq p/\sqrt{(k,p-1)}$ (q.v. \cite{Ak1}) -- although the map $\eta_{a,k}$ is not a permutation whenever this result is nontrivial!  Surprisingly, extensive computer evidence generated by the author strongly suggests a much better result than Theorem \ref{almostpower}, which we consider our most intriguing conjecture.

\begin{conj} For all $k \in \mathbb{Z}_p^\times$ with $(k,p-1)=1$ and $k \geq 2$,
$$
S(a,k,M) \ll p^{3/4},
$$
uniformly in $a$, $k$, and $M$.
\end{conj}

All the standard techniques appear not to help at all with this question.

\subsection{S\'{o}s Permutations} \label{Irrational}

Recall the definition of the S\'{o}s permutation $\beta_\alpha$:

\begin{definition} For $s,t \in \{1,\ldots,n\}=[n]$ and $\alpha \in \mathbb{R}$ irrational, $\beta_\alpha(s) < \beta_\alpha(t)$ iff $\{\alpha s\} < \{\alpha t\}$, where $\{x\}$ is the fractional part of $x$.
\end{definition}

Equivalently, $\beta_\alpha(t)$ is the number of $s \in [n]$ with $\{\alpha s\} \leq \{\alpha t\}$.  It is clear that
$$
D^*(\beta_\alpha) = \max_{s,t \in [n]} \left | \left | \beta_\alpha([s]) \cap [t] \right | - \frac{st}{n} \right |.
$$
The cardinality of the set $\beta_\alpha([s]) \cap [t]$ is the number of $x$'s in $[s]$ so that the number of $y$'s in $[n]$ with $\{\alpha y\} \leq \{\alpha x\}$ is less than or equal to $t$.  If we let $\mathcal{A}_s(\alpha) = \{\{\alpha x\} : x \in [s]\}$, then we wish to know the maximum value of
$$
\left|\left|\mathcal{A}_s(\alpha) \cap [0,\{\alpha \beta_\alpha^{-1}(t)\}]\right|-\frac{st}{n}\right|
$$
over all $s,t \in [0,n-1]$, since $\{\alpha \beta_\alpha^{-1}(t)\}$ is the $t^\text{th}$ smallest point of $\mathcal{A}_n(\alpha)$.  Now, for a set of reals $A \subset [0,1]$, define
$$
d^*(A) = \sup_{0 \leq x \leq 1} \big||A \cap [0,x]|-x|A|\big|.
$$
We can now write
\begin{align*}
D^*(\beta_\alpha) & = \max_{s,t \in [n]} \left|\left|\mathcal{A}_s(\alpha) \cap [0,\{\alpha \beta_\alpha^{-1}(t)\}]\right|-\frac{st}{n}\right| \\
& = \max_{s,t \in [n]} \left|\left|\mathcal{A}_s(\alpha) \cap [0,\{\alpha \beta_\alpha^{-1}(t)\}]\right|-s\{\alpha \beta_\alpha^{-1}(t)\} + s\{\alpha \beta_\alpha^{-1}(t)\} - \frac{st}{n}\right| \\
& \leq \max_{s \in [n]} d^*(\mathcal{A}_s(\alpha)) + s \max_t \left | \{\alpha \beta_\alpha^{-1}(t)\} - \frac{t}{n} \right | \\
& = \max_{s \in [n]} d^*(\mathcal{A}_s(\alpha)) + \frac{s}{n} \max_t \left | n\{\alpha t\} - \beta_\alpha(t) \right |
\end{align*}
But, since $\beta_\alpha(t) = |\mathcal{A}_s(\alpha) \cap [0,\{\alpha t\}]|$, we have
\begin{equation} \label{discrelation}
D^*(\beta_\alpha) \leq \max_{s \in [n]} d^*(\mathcal{A}_s(\alpha)) + d^*(\mathcal{A}_n(\alpha)) \leq 2 \max_{s \in [n]} d^*(\mathcal{A}_s(\alpha)).
\end{equation}

An old and well-known result of H. Weyl states that $d^*(\mathcal{A}_n(\alpha))$ is $o(n)$ for any $\alpha$ irrational.  Therefore, we may conclude

\begin{theorem} If $\alpha$ is irrational, then $\beta_\alpha$ is quasirandom.
\end{theorem}

It is a theorem of J. Schoi{\ss}engeier (\cite{Sg1}) that $d^*(\mathcal{A}_n(\alpha)) \ll \log n$ iff the partial quotients of the continued fraction of $\alpha$ are bounded in average.  Therefore we have

\begin{theorem} If $\alpha$ is irrational and has partial quotients bounded in average, then $\beta_\alpha$ is maximally quasirandom.
\end{theorem}

Furthermore, it is a theorem of Khintchine (q.v. \cite{DT1}) that
$$
\max_{1 \leq s \leq n} d^*(\mathcal{A}_s(\alpha)) \ll \log n \cdot f(\log n)
$$
for almost all $\alpha$ if and only if
$$
\sum_{n=1}^\infty \frac{1}{n f(n)} < \infty.
$$
So we also have
\begin{cor} If $\sum_{n=1}^\infty (n f(n))^{-1} < \infty$, then $D(\beta_\alpha) \ll \log n \cdot f(\log n)$ for almost all $\alpha \in \mathbb{R}$.
\end{cor}

Therefore, it is clear that $D(\beta_\alpha) \ll \log n (\log \log n)^{1+\epsilon}$ almost always.  We return to Conjecture \ref{weaklogn} now.  Define the {\it continuant} $K(a_1,a_2,\ldots,a_m)$ to be the denominator of the continued fraction $p/q = [0;a_1,a_2,\ldots,a_m]$, and define $\mathcal{F}(B)$, for each $B \geq 1$, to be the set of continuants of sequences of partial quotients bounded in average by $B$.

\begin{prop} \label{biaimplieslog} For any $n \in \mathcal{F}(B)$,
$$
\min_{k \in [n]} D(\psi_k) \ll \log n.
$$
where the implicit constant depends only on $B$.
\end{prop}

\begin{proof} It is implicit in the work of Schoi{\ss}engeier (\cite{Sg1}) (and explicit in \cite{DT1}) that there exists absolute constants $C$ and $N$ so that if $n>N$, and the irrational $\alpha$ has continued fraction expansion $[a_0;a_1,a_2,\ldots]$ with convergents $\{p_s/q_s\}_{s\geq1}$, then
$$
d^*(\{\alpha s\}_{1\leq s\leq n}) \leq C (\sum_{i=1}^{m} a_i + m)
$$
where $m$ is chosen so that $q_m \leq n \leq q_{m+1}$.  Since $m \ll \log n$, the right hand side is $\leq (CB+1)m \ll \log n$.
\end{proof}

Zaremba's Conjecture (\cite{Z}) implies that $\mathcal{F}_5 = \mathbb{N}$ (!) -- so this clearly implies Conjecture \ref{weaklogn} if it is true.  We can ask for considerably less, however:

\begin{prop} Choose $B \geq 2$, and let $\mathcal{C}_B$ be the set of irrationals whose partial quotients are bounded in average by $B$.  If
\begin{equation} \label{logbnd}
\inf_{\alpha \in \mathcal{C}_B} \min_{k \in [n]} \left | \frac{k}{n} - \alpha \right | \ll \frac{\log n}{n^2}
\end{equation}
then Conjecture \ref{weaklogn} follows.
\end{prop}
\begin{proof} Fix $n$ sufficiently large.  If (\ref{logbnd}) holds, we may choose $k \in [n]$ and $\alpha \in \mathcal{C}_B$ so that $|k/n - \alpha| \leq C n^{-2} \log n$ for some $C > 0$.  By (\ref{discrelation}) it suffices to prove that $d^*(\mathcal{A}_s(k/n)) \ll \log n$.  Define $I_t = [0,t/n]$.  Then, we have
\begin{align*}
d^*(\mathcal{A}_s(k/n)) &= \max_{t\in [n]} \left | \left | \mathcal{A}_s(k/n) \cap I_t \right | - \frac{st}{n} \right | \\
& = \max_{t\in [n]} \left | \left | \mathcal{A}_s(k/n) \cap I_t \right | - \left | \mathcal{A}_s(\alpha) \cap I_t \right | + \left | \mathcal{A}_s(\alpha) \cap I_t \right | - \frac{st}{n} \right | \\
& \leq \max_{t\in [n]} \big | \left | \mathcal{A}_s(k/n) \cap I_t \right | - \left | \mathcal{A}_s(\alpha) \cap I_t \right | \big | + \max_{t\in [n]} \left | \left | \mathcal{A}_s(\alpha) \cap I_t \right | - \frac{st}{n} \right | \\
& \leq \max_{t\in [n]} \big | \left | \mathcal{A}_s(k/n) \cap I_t \right | - \left | \mathcal{A}_s(\alpha) \cap I_t \right | \big | + d^*(\mathcal{A}_s(\alpha))
\end{align*}
The second summand is $\ll \log n$ by Proposition \ref{biaimplieslog}.  The first summand, for a given $t$, is the number of multiples (up to $s$) of $k/n$ lying in $I_t$ minus the number of multiples (up to $s$) of $\alpha$ lying in $I_t$.  If $\{kj/n\} \in [\epsilon,t-\epsilon]$, then $\{\alpha j\} \in I_t$, if $t > \epsilon = |k/n-\alpha|\cdot n$.  But, $|k/n-\alpha|\cdot n \leq C n^{-1} \log n$, so this quantity is bounded by the number $N_1$ of points $\{\alpha j\} \in [0,\epsilon) \cup (t-\epsilon,t]$ plus the number $N_2$ of points $\{kj/n\} \in [0,\epsilon) \cup (t-\epsilon,t]$, where $j$ varies from $1$ to $n$.  It is easy to see that Proposition \ref{biaimplieslog} gives $N_1 \ll \log n$ and the fact that $N_2 \ll \log n$ is trivial.
\end{proof}

We cannot prove that (\ref{logbnd}) holds, although we believe it to be true for some very small $B$.  Therefore, we have the following conjecture.

\begin{conj} For some integer $B \geq 2$,
$$
\inf_{\alpha \in \mathcal{C}_B} \min_{k \in [n]} \left | \frac{k}{n} - \alpha \right | \ll \frac{\log n}{n^2}
$$
\end{conj}

\section{Application: Simultaneous Initial Intervals of S\'{o}s Permutations} \label{appsfrac}

In \cite{My1}, K. O'Bryant asks, for a given irrational $\alpha$, what values can be taken on by the following function:
$$
B_\alpha(k) = \left | \{ 1 \leq q \leq k : \{q\alpha\} \leq \{k \alpha\} \} \right |.
$$
Call the set of all such cardinalities $A_\alpha$.  In particular, O'Bryant poses a series of questions:

\begin{enumerate}
\item Can $A_\alpha$ ever be $\mathbb{Z}^+$?
\item Is $8 \not \in A_{\sqrt{2}}$?  Is $4 \not \in A_{-(1+\sqrt{5})/2}$?
\item Is it true that if the continued fraction expansion of $\alpha$ has bounded partial quotients, then $A_\alpha$ has positive density in the naturals?
\end{enumerate}

The answer to the first question, as noted in \cite{My1}, is actually ``yes''.  This follows by choosing any $\alpha$ whose odd-numbered partial quotients are unbounded.  We address the third question here using the theory of quasirandom permutations.  First, generalize $B_\alpha(k)$ and $A_\alpha$ as follows: for a permutation $\sigma$ of $[n]$ and $k \in [n]$, define
$$
B_\sigma(k) = \left | \{ 1 \leq q \leq k : \sigma(q) \leq \sigma(k) \right |
$$
and
$$
A_\sigma = \{B_\sigma(k) : k \in [n]\}
$$
Clearly, then, the definitions of $A_{\beta_\alpha}$ and $A_\alpha$ agree and definitions of $B_{\beta_\alpha}(k)$ and $B_\alpha(k)$ agree, so long as $n \geq k$.  We also define $D_0(\sigma)$ as
$$
D_0(\sigma) = \max_{I,J} D_J(\sigma(I)),
$$
where $I$ and $J$ are intervals of $[n]$ (i.e., they do not ``wrap around'').  We have the following general result.

\begin{theorem} Given $\sigma$ a permutation of $[n]$, there is an element of $A_\sigma$ in every interval of length $\sqrt{32nD(\sigma)}$ lying in $[n]$.
\end{theorem}
\begin{proof} Given $r \leq n - 4\sqrt{nD_0(\sigma)}$, we wish to know if there exists a $k$ so that
$$
\left | \{ 1 \leq q \leq k : \sigma(q) \leq \sigma(k) \} \right |
$$
is at least $r$ but no larger than $r + 4\sqrt{nD_0(\sigma)}$.  The question is equivalent to asking whether there exists a point in the interval $[r,r+4\sqrt{nD_0(\sigma)}]$ equal to
$$
B_\sigma(k) = \left | \sigma([k]) \cap [\sigma(k)] \right |
$$
for some $k$.  By the definition of discrepancy,
\begin{equation} \label{BBound}
\left | B_\sigma(k) - \frac{k \sigma(k)}{n} \right | \leq D_0(\sigma).
\end{equation}
Define $r^\prime = r + 3\sqrt{nD_0(\sigma)}$, let $S=\lfloor \sqrt{nD_0(\sigma)} \rfloor+1$ and $t = \lfloor \sqrt{r^\prime n} \rfloor$, and define $I = [t,t+S-1]$.  So long as $t + S - 1\leq n$, which is true since $r^\prime \leq n - \sqrt{nD_0(\sigma)}$, this is a genuine subset of $[n]$, and so
$$
\left | |\sigma(I) \cap I| - \frac{S^2}{n} \right | \leq D_0(\sigma).
$$
Therefore, there is some $k \in I$ whose image under $\sigma$ lies in $I$, since $S^2/n > D_0(\sigma)$.  We have
\begin{align*}
\left |\frac{k \sigma(k)}{n} - r^\prime \right | &\leq \frac{\left (\lfloor \sqrt{r^\prime n}\rfloor + \lfloor \sqrt{nD_0(\sigma)}\rfloor\right)^2}{n}- r^\prime \\
& \leq 2 \sqrt{r^\prime D(\sigma)} + D(\sigma) \\
& \leq 3 \sqrt{n D_0(\sigma)}.
\end{align*}
Combining this with (\ref{BBound}), we have $B_\sigma(k) \in [r,r+4 \sqrt{nD_0(\sigma)}]$.  Since $D(\sigma) \leq 2D_0(\sigma)$, the result follows.
\end{proof}

We bring the discussion back to S\'{o}s permutations.

\begin{cor} If $\alpha$ has partial quotients bounded in average, then $|A_\alpha \cap [n]| \gg \sqrt{n / \log n}$.
\end{cor}
\begin{proof} The statement follows immediately from Proposition \ref{biaimplieslog} and the previous theorem.
\end{proof}

\section*{Acknowledgements}
Thank you to Fan Chung and Ron Graham for their endless support, and for suggesting the problems which eventually motivated the current work.  I also wish to thank Kevin O'Bryant for introducing me (thoroughly!) to S\'{o}s permutations.  Gratitude also to J\'{o}zsef Beck, Ron Evans, Michael Lacey, Heather Stoneberg-Henry, and Lei Wu for helpful discussions.

\end{document}